\documentclass[fleqn]{mat01}
\usepackage{times,mathtimy,amssymb,latexsym,amsbsy}
\begin{document}

\setcounter{page}{117}
\firstpage{117}

\newtheorem{theore}{Theorem}
\renewcommand\thetheore{\arabic{section}.\arabic{theore}}
\newtheorem{theor}{\bf Theorem}
\newtheorem{rem}[theore]{Remark}
\newtheorem{propo}[theore]{\rm PROPOSITION}
\newtheorem{lem}[theore]{Lemma}
\newtheorem{definit}[theore]{\rm DEFINITION}
\newtheorem{coro}{\rm COROLLARY}
\newtheorem{exampl}[theore]{Example}
\newtheorem{case}{Case}

\def\corol{\trivlist \item[\hskip \labelsep{COROLLARY.}]}
\def\noteproof{\trivlist \item[\hskip \labelsep{\it Note added in Proof.}]}

\title{Irreducibility results for compositions of
polynomials in~several variables}

\markboth{Anca Iuliana Bonciocat and Alexandru
Zaharescu}{Compositions of polynomials in several variables}

\author{ANCA IULIANA BONCIOCAT$^{1,2}$ and ALEXANDRU~ZAHARESCU$^{3}$}

\address{$^{1}$Universit\"{a}t Bonn, Institut f\"{u}r Angewandte Mathematik,
Abt. f\"{u}r Stochastik, 53115~Bonn, Wegelerstr. 6, Deutschland\\
\noindent $^{2}$Institute of Mathematics ``Simion Stoilow'' of the
Romanian Academy, P.O.~Box~1-764, RO-70700 \ Bucharest,
Romania\\
\noindent $^{3}$Department of Mathematics, University of
Illinois at Urbana-Champaign, Altgeld~Hall, 1409 W. Green Street,
Urbana, IL 61801, USA\\
\noindent E-mail: anca@wiener.iam.uni-bonn.de;
Anca.Bonciocat@imar.ro; zaharesc@math.uiuc.edu}

\volume{115}

\mon{May}

\parts{2}

\pubyear{2005}

\Date{MS received 28 September 2004; revised 15 January 2005}

\begin{abstract}
We obtain explicit upper bounds for the number of irreducible
factors for a class of compositions of polynomials in several
variables over a given field. In particular, some irreducibility
criteria are given for this class of compositions of polynomials.
\end{abstract}

\keyword{Composition of polynomials; irreducibility results.}

\maketitle

\section{Introduction}

In connection with Hilbert's irreducibility theorem, Cavachi
proved in \cite{Cav} that for any relatively prime polynomials
$f(X),g(X)\in {\mathbb{Q}}[X]$ with $\deg f<\deg g$, the
polynomial $f(X)+pg(X)$ is irreducible over ${\mathbb{Q}}$ for all
but finitely many prime numbers $p$. Sharp explicit upper bounds
for the number of factors over ${\mathbb{Q}}$ of a linear
combination $n_{1}f(X)+n_{2}g(X)$, covering also the case $\deg f
= \deg g$, have been derived in \cite{Boncio}. In \cite{BonZah},
we realized that by using technics similar to those employed in
\cite{CVZ} and \cite{Boncio}, upper bounds for the number of
factors and irreducibility results can also be obtained for a
class of compositions of polynomials of one variable with integer
coefficients. More specifically, the following result is proved in
\cite{BonZah}.\vspace{.7pc}

\noindent Let $f(X)=a_{0}+a_{1}X+\cdots +a_{m}X^{m}$ and
$g(X)=b_{0}+b_{1}X+\cdots +b_{n}X^{n}\in {\mathbb{Z}}[X]$ be
nonconstant polynomials of degree $m$ and $n$ respectively, with
$a_{0}\neq 0$, and let $L_{1}(f)=|a_{0}|+\cdots +|a_{m-1}|$.
Assume that $d_{1}$ is a positive divisor of $a_{m}$ and $d_{2}$ a
positive divisor of $b_{n}$ such that
\begin{equation*}
|a_{m}|>d_{1}^{mn}d_{2}^{m^{2}n}L_{1}(f).
\end{equation*}
Then the polynomial $f\circ g$ has at most $\Omega
(a_{m}/d_{1})+m\Omega (b_{n}/d_{2})$ irreducible factors over
${\mathbb{Q}}$, where $\Omega (k)$ is the total number of prime
factors of $k$, counting multiplicities. The same conclusion holds
in the wider range
\begin{equation*}
|a_{m}|>d_{1}^{n}d_{2}^{mn}L_{1}(f),
\end{equation*}
provided that $f$ is irreducible over $\mathbb{Q}$.\newpage

In the present paper we provide explicit upper bounds for the
number of factors, and irreducibility results for a class of
compositions of polynomials in several variables over a given
field. We will deduce this result from the corresponding result
for polynomials in two variables $X,Y$ over a field $K$. We use
the following notation. For any polynomial $f\in K[X,Y]$ we denote
by $\deg _{Y}f$ the degree of $f$ as a polynomial in $Y$, with
coefficients in $K[X]$. Then we write any polynomial $f\in K[X,Y]$
in the form
\begin{equation*}
f=a_{0}(X)+a_{1}(X)Y+\cdots +a_{d}(X)Y^{d},
\end{equation*}
with $a_{0},a_{1},\ldots ,a_{d}$ in $K[X]$, $a_{d}\neq 0$, and define
\begin{equation*}
H_{1}(f)=\max \{\deg a_{0},\ldots ,\deg a_{d-1}\}.
\end{equation*}
Finally, for any polynomial $f\in K[X]$ we denote by $\Omega (f)$
the number of irreducible factors of $f$, counting multiplicities
($\Omega (c)=0$ for $c\in K$). We will prove the following
theorem.

\begin{theor}[\!]
Let $K$ be a field and let $f(X,Y)=a_{0}+a_{1}Y+\cdots
+a_{m}Y^{m}${\rm ,} $g(X,Y)=b_{0}+b_{1}Y+\cdots +b_{n}Y^{n}${\rm
,} with $a_{0},a_{1},\ldots ,a_{m}${\rm ,} $b_{0},b_{1},\ldots
,b_{n}\in K[X]${\rm ,} $a_{0}a_{m}b_{n}\neq 0$. If $d_{1}$ is a
factor of $a_{m}$ and $d_{2}$ a factor of $b_{n}$ such that
\begin{equation*} \deg
a_{m}>mn\deg d_{1}+m^{2}n\deg d_{2}+H_{1}(f),
\end{equation*}
then the polynomial $f(X,g(X,Y))$ has at most $\Omega
(a_{m}/d_{1})+m\Omega (b_{n}/d_{2})$ irreducible factors over
$K(X)$. The same conclusion holds in the wider range
\begin{equation*}
\deg a_{m}>n\deg d_{1}+mn\deg d_{2}+H_{1}(f),
\end{equation*}
provided that $f$ is irreducible over $K(X)$.
\end{theor}
Theorem~1 provides, in particular, bounds for the number of
irreducible factors of $f(X,Y)$ over $K(X)$, by taking $g(X,Y)=Y$.

\begin{coro}$\left.\right.$\vspace{.5pc}

\noindent {\it Let $K$ be a field and let
$f(X,Y)=a_{0}+a_{1}Y+\cdots + a_{m}Y^{m},$ with
$a_{0},a_{1},\ldots, a_{m}\in K[X]${\rm ,} $a_{0}a_{m}\neq 0$. If
$d$ is a factor of $a_{m}$ such that
\begin{equation*}
\deg a_{m}>m\deg d+H_{1}(f),
\end{equation*}
then the polynomial $f(X,Y)$ has at most $\Omega (a_{m}/d)$
irreducible factors over $K(X)$.}
\end{coro}
Under the assumption that $a_{m}$ has an irreducible factor over
$K$ of large enough degree, we have the following irreducibility
criteria.

\begin{coro}$\left.\right.$\vspace{.5pc}

\noindent {\it Let $K$ be a field and let
$f(X,Y)=a_{0}+a_{1}Y+\cdots +a_{m}Y^{m}${\rm ,} with
$a_{0},a_{1},\ldots,a_{m}\in K[X]${\rm ,} $a_{0}a_{m}\neq 0$. If
$a_{m}=pq$ with $p,q\in K[X]${\rm ,} $p$ irreducible over $K${\rm
,} and
\begin{equation*}
\deg p>(m-1)\deg q+H_{1}(f),
\end{equation*}
then the polynomial $f(X,Y)$ is irreducible over $K(X)$.}
\end{coro}

\begin{coro}$\left.\right.$\vspace{.5pc}

\noindent {\it Let $K$ be a field and let
$f(X,Y)=a_{0}+a_{1}Y+\cdots +a_{m}Y^{m}${\rm ,}
$g(X,Y)=b_{0}+b_{1}Y+\cdots +b_{n}Y^{n}${\rm ,} with
$a_{0},a_{1},\ldots,a_{m}, b_{0},b_{1},\ldots ,b_{n}\in K[X]${\rm
,} $a_{0}a_{m}b_{n}\neq 0${\rm ,} and $f$ irreducible over $K(X)$. If
$a_{m}=pq$ with $p,q\in K[X], p$ irreducible over $K${\rm ,} and
\begin{equation*}
\deg p>(n-1)\deg q+mn\deg b_{n}+H_{1}(f),
\end{equation*}
then the polynomial $f(X,g(X,Y))$ is irreducible over $K(X)$.}
\end{coro}\vspace{.3pc}

\begin{coro}$\left.\right.$\vspace{.5pc}

\noindent {\it Let $K$ be a field and let
$f(X,Y)=a_{0}+a_{1}Y+\cdots +a_{m}Y^{m}${\rm ,}
$g(X,Y)=b_{0}+b_{1}Y+\cdots +b_{n}Y^{n}${\rm ,} with
$a_{0},a_{1},\ldots ,a_{m}, b_{0},b_{1},\ldots ,b_{n}\in K[X]${\rm ,}
$a_{0}a_{m}b_{n}\neq 0$. If
$a_{m}=pq$ with $p,q\in K[X], p$ irreducible over $K${\rm ,} and
\begin{equation*}
\deg p>\max \{(m-1)\deg q,\ (n-1)\deg q+mn\deg b_{n}\}+H_{1}(f),
\end{equation*}
then the polynomial $f(X,g(X,Y))$ is irreducible over $K(X)$.}\medskip
\end{coro}
Another consequence of Theorem~1 is the following corresponding
result for polynomials in $r\geq 2$ variables $X_{1},X_{2},\ldots
,X_{r}$ over $K$. In this case, for any polynomial $f\in
K[X_{1},\ldots ,X_{r}]$, $\Omega (f)$ will stand for the number of
irreducible factors of $f$ over $K(X_{1},\ldots ,X_{r-1})$,
counting multiplicities. Then, for any polynomial $f\in
K[X_{1},\ldots ,X_{r}]$ and any $j\in \{1,\ldots ,r\}$ we denote
by $\deg _{X_{j}}f$ the degree of $f$ as a polynomial in $X_{j}$
with coefficients in $K[X_{1},\ldots ,X_{j-1},X_{j+1},\ldots
,X_{r}]$. We also write any polynomial $f\in K[X_{1},\ldots
,X_{r}]$ in the form
\begin{align*}
f&=a_{0}(X_{1},\ldots ,X_{r-1})+a_{1}(X_{1},\ldots
,X_{r-1})X_{r}\\[.3pc]
&\quad\,+\cdots + a_{d}(X_{1},\ldots ,X_{r-1})X_{r}^{d},
\end{align*}
with $a_{0},a_{1},\ldots ,a_{d}\in K[X_{1},\ldots ,X_{r-1}]$,
$a_{0}\neq 0$, and for any $j\in \{1,\ldots ,r-1\}$ we let
\begin{equation*}
H_{j}(f)=\max \{\deg _{X_{j}}a_{0},\deg _{X_{j}}a_{1},\ldots ,\deg
_{X_{j}}a_{d-1}\}.
\end{equation*}
Then one has the following result.\vspace{.5pc}

\begin{coro}$\left.\right.$\vspace{.5pc}

\noindent {\it Let $K$ be a field{\rm ,} $r\geq 2${\rm ,} and let
$f(X_{1},\ldots ,X_{r})=a_{0}+a_{1}X_{r}+\cdots
+a_{m}X_{r}^{m}${\rm ,} $g(X_{1},\ldots
,X_{r})=b_{0}+b_{1}X_{r}+\cdots +b_{n}X_{r}^{n}${\rm ,} with
$a_{0},a_{1},\ldots,a_{m}${\rm ,} $b_{0},b_{1},\ldots ,b_{n}\in
K[X_{1},\ldots ,X_{r-1}]${\rm ,} $a_{0}a_{m}b_{n}\neq 0$. If $d_{1}$ is a factor of $a_{m}$ and $d_{2}$ a
factor of $b_{n}$ such that for an index $j\in \{1,\ldots
,r-1\}${\rm ,}
\begin{equation*}
\deg _{X_{j}}a_{m}>mn\deg _{X_{j}}d_{1}+m^{2}n\deg _{X_{j}}d_{2}+H_{j}(f),
\end{equation*}
then the polynomial $f(X_{1},\ldots ,X_{r-1},g(X_{1},\ldots
,X_{r}))$ has at most $\Omega (a_{m}/d_{1})+m\Omega
(b_{n}/d_{2})$\ irreducible factors over the field $K(X_{1},\ldots
,X_{r-1})$. The same conclusion holds in the wider range
\begin{equation*}
\deg _{X_{j}}a_{m}>n\deg _{X_{j}}d_{1}+mn\deg _{X_{j}}d_{2}+H_{j}(f),
\end{equation*}
provided that $f$ is irreducible over $K(X_{1},\ldots ,X_{r-1})$.}
\end{coro}
In particular we have the following irreducibility
criterion.\vspace{.3pc}

\begin{coro}$\left.\right.$\vspace{.5pc}

\noindent {\it Let $K$ be a field{\rm ,} $r\geq 2${\rm ,} and let
$f(X_{1},\ldots ,X_{r})=a_{0}+a_{1}X_{r}+\cdots
+a_{m}X_{r}^{m}${\rm ,} $g(X_{1},\ldots
,X_{r})=b_{0}+b_{1}X_{r}+\cdots +b_{n}X_{r}^{n}${\rm ,} with
$a_{0},a_{1},\ldots ,a_{m}, b_{0},b_{1},\ldots ,b_{n}\in
K[X_{1},\ldots ,X_{r-1}], a_{0}a_{m}b_{n}\neq 0$. If $a_{m}=p\cdot q$ with $p$ a prime element of
the ring $K[X_{1},\ldots ,X_{r-1}]$ such that for an index $j\in
\{1,\ldots ,r-1\}${\rm ,}
\begin{align*}
\deg _{X_{j}}p>\max \{(m-1)\deg _{X_{j}}q,(n-1) \deg
_{X_{j}}q+mn\deg _{X_{j}}b_{n}\}\!+\!H_{j}(f),
\end{align*}
then the polynomial $f(X_{1},\ldots ,X_{r-1},g(X_{1},\ldots
,X_{r}))$ is irreducible over the field $K(X_{1},\ldots
,X_{r-1})$.}
\end{coro}
Corollary 5 follows from Theorem~1 by writing $Y$ for $X_{r}$ and
$X$ for $X_{j}$, where $j$ is any index for which
\begin{equation*}
\deg _{X_{j}}a_{m}>mn\deg _{X_{j}}d_{1}+m^{2}n\deg _{X_{j}}d_{2}+H_{j}(f),
\end{equation*}
and by replacing the field $K$ with the field generated by $K$ and
the variables $X_{1},X_{2},\ldots,$ $X_{r-1}$ except for $X_{j}$.

The reader may naturally wonder how sharp the above results are.
In this connection, we discuss a couple of examples in the next
section.

\section{Examples}

Let $K={\mathbb{Q}}$, choose integers $m,d\geq 2$, select
polynomials $a_{0}(X),a_{1}(X),\ldots ,a_{m-1}(X)\in
{\mathbb{Q}}[X]$ with $a_{0}(X)\neq 0$, and consider the
polynomial in two variables $f(X,Y)$ given by
\begin{align*}
f(X,Y)=a_{0}(X)+a_{1}(X)Y+\cdots
+a_{m-1}(X)Y^{m-1}+(X^{d}+5X+5)Y^{m}.
\end{align*}
Under these circumstances, in terms of the degrees of the
polynomials $a_{0}(X), a_{1}(X),$ $\ldots ,a_{m-1}(X)$, can we be
sure that the polynomial $f(X,Y)$ is irreducible over
${\mathbb{Q}}(X)$? The polynomial $p(X)=X^{d}+5X+5$ is an
Eisensteinian polynomial with respect to the prime number $5$, and
hence it is irreducible over ${\mathbb{Q}}$. We may then apply
Corollary 2, with $q=1$, in order to conclude that $f(X,Y)$ is
irreducible over ${\mathbb{Q}}(X)$ as long as $H_{1}(f)<d$, that
is, as long as each of the polynomials $a_{0}(X),$ $a_{1}(X),$
$\ldots,a_{m-1}(X)$ has degree less than or equal to $d-1$. We
remark that for any choice of $m,d\geq 2$ this bound is the best
possible, in the sense that there are polynomials
$a_{0}(X),a_{1}(X),\ldots ,a_{m-1}(X)\in {\mathbb{Q}}[X]$,
$a_{0}(X)\neq 0$, for which
\begin{equation*}
\max \{\deg _{X}a_{0}(X),\deg _{X}a_{1}(X),\ldots ,\deg _{X}a_{m-1}(X)\}=d,
\end{equation*}
such that the corresponding polynomial $f(X,Y)$ is reducible over
${\mathbb{Q}}(X)$. Indeed, one may choose for instance
$a_{0}(X),a_{1}(X),\ldots ,a_{m-2}(X)$ to be any polynomials with
coefficients in ${\mathbb{Q}}$, with $a_{0}(X)\neq 0$, of degrees
less than or equal to $d$, and define $a_{m-1}(X)$ by the equality
\begin{equation*}
a_{m-1}(X)=-X^{d}-5X-5-\sum_{0\leq i\leq m-2}a_{i}(X).
\end{equation*}
Then, on the one hand, we will have $\max \{\deg _{X}a_{0},\ldots
,\deg _{X}a_{m-1}\}=d$ and on the other hand, the corresponding
polynomial $f(X,Y)$ will be reducible over ${\mathbb{Q}}(X)$,
being divisible by $Y-1$.

In the next example, let us slightly modify the polynomial
$f(X,Y)$, and choose a polynomial $g(X,Y)$ of arbitrary degree,
say
\begin{align*}
f(X,Y) &= a_{0}+a_{1}Y+\cdots +a_{m-1}Y^{m-1}
+(X^{n}+5X+5)^{2}Y^{m},\\[.3pc]
g(X,Y) &= b_{0}+b_{1}Y+\cdots +b_{n-1}Y^{n-1}+Y^{n},
\end{align*}
where $a_{0},a_{1},\ldots ,a_{m-1},\ b_{0},b_{1},\ldots
,b_{n-1}\in {\mathbb{Q}}[X]$, $a_{0}(X)\neq 0$. We may apply
Theorem~1, with $d_{1}=d_{2}=1$, in order to conclude that
$f(X,g(X,Y))$ has at most two irreducible factors over
${\mathbb{Q}}(X)$ as long as $H_{1}(f)<2n$, that is, as long as
each of the polynomials $a_{0},a_{1},\ldots ,a_{m-1}$ has degree
less than or equal to $2n-1$. This bound too is the best possible,
as there exist polynomials  $a_{0},a_{1},\ldots ,a_{m-1}\in
{\mathbb{Q}}[X]$, $a_{0}(X)\neq 0$, $g\in {\mathbb{Q}}[Y]$, $g$
monic, for which
\begin{equation*}
\max \{\deg a_{0},\deg a_{1},\ldots ,\deg a_{m-1}\}=2n,
\end{equation*}
such that the corresponding polynomial $f(X,g(X,Y))$ has at least
three irreducible factors over ${\mathbb{Q}}(X)$. For instance,
one may take $g(Y)=Y^{2}$, choose polynomials
$a_{0}(X),a_{1}(X),\ldots ,a_{m-2}(X)$ with coefficients in
${\mathbb{Q}}$, with $a_{0}(X)\neq 0$, of degrees less than or
equal to $2n$, and define $a_{m-1}(X)$ by the equality
\begin{equation*}
a_{m-1}(X)=-(X^{n}+5X+5)^{2}-\sum_{0\leq i\leq m-2}a_{i}(X).
\end{equation*}
Then we will have $\max \{\deg a_{0},\ldots ,\deg a_{m-1}\}=2n$,
while the corresponding polynomial $f(X,g(X,Y))$ will have at
least three irreducible factors over ${\mathbb{Q}}(X)$, being
divisible by $Y^{2}-1$.

\section{Proof of Theorem 1}

Let $K$, $f(X,Y)$, $g(X,Y)$, $d_{1}$ and $ d_{2}$ be as
in the statement of the theorem. Let $b\in K[X]$ denote
the greatest common divisor of $b_{0},b_{1},\ldots
,b_{n}$, and define the polynomial $\overline{g}(X,Y)\in
K[X,Y]$ by the equality
\begin{equation*}
g(X,Y)=b_{0}+b_{1}Y+\cdots +b_{n}Y^{n}=b\overline{g}(X,Y).
\end{equation*}
Next, let $a\in K[X]$ denote the greatest common divisor of the
coefficients of $f(X,g(X,Y))$ viewed as a polynomial in $Y$, and
define the polynomial $ F(X,Y)\in K[X,Y]$ by the equality
\begin{equation*}
f(X,g(X,Y))=aF(X,Y).
\end{equation*}
If we assume that $f(X,g(X,Y))$ has $s>\Omega
(a_{m}/d_{1})+m\Omega (b_{n}/d_{2})$ irreducible factors over
$K(X)$, then the polynomial $F(X,Y)$ will have a factorization
$F(X,Y)=F_{1}(X,Y)\cdots F_{s}(X,Y)$, with $F_{1}(X,Y), \ldots,
F_{s}(X,Y)\in K[X,Y], \deg _{Y}F_{1}(X,Y)\geq 1, \ldots , \deg
_{Y}F_{s}(X,Y)\geq 1$. Let $t_{1}, \ldots , t_{s}\in K[X]$ be the
leading coefficients of $F_{1}(X,Y), \ldots , F_{s}(X,Y)$
\,respectively, \,viewed \,as \,polynomials \,in \,$Y$. \,By
\,comparing the

\noindent leading coefficients in the equality
\begin{equation*}
a_{0}+a_{1}g(X,Y)+\cdots +a_{m}g^{m}(X,Y)=aF_{1}(X,Y)\cdots F_{s}(X,Y)
\end{equation*}
we obtain the following equality in $K[X]$:
\begin{equation}
at_{1}\cdots t_{s}=a_{m}b_{n}^{m}=d_{1}d_{2}^{m}\cdot
\frac{a_{m}}{d_{1}} \cdot \left( \frac{b_{n}}{d_{2}}\right) ^{m}.
\label{ec1}
\end{equation}
Then, in view of (\ref{ec1}), it follows easily that at least one
of the $t_{i}$'s, say $t_{1}$, will divide $d_{1}d_{2}^{m}$. As a
consequence, one has
\begin{equation}
\deg t_{1}\leq \deg d_{1}+m\deg d_{2}.  \label{ec2}
\end{equation}
We now consider the polynomial
\begin{align*}
h(X,Y) &= f(X,g(X,Y))-a_{m}(X)g(X,Y)^{m} \\[.3pc]
&= a_{0}(X)+a_{1}(X)g(X,Y)+\cdots +a_{m-1}(X)g(X,Y)^{m-1}.
\end{align*}
Recall that $a_{0}$ and $\overline{g}$ are relatively prime, and
$a_{0}(X)\neq 0$. It follows that the polynomials
$\overline{g}^{m}(X,Y)$ and $h(X,Y)$ are relatively prime.
Therefore $\overline{g}^{m}(X,Y)$ and $F_{1}(X,Y)$ are relatively
prime. As a consequence, the resultant $R(\overline{g}^{m},F_{1})$
of $\overline{g}^{m}(X,Y)$ and $F_{1}(X,Y)$, viewed as polynomials
in $Y$ with coefficients in $K[X]$, will be a nonzero element of
$K[X]$. We now introduce a nonarchimedean absolute value $|\cdot
|$ on $K(X)$, as follows. We fix a real number $\rho $, with
$0<\rho <1$, and for any polynomial $F(X)\in K[X]$ we define
$|F(X)|$ by the equality
\begin{equation}
|F(X)|=\rho ^{-\deg F(X)}.  \label{ec3}
\end{equation}
We then extend the absolute value $|\cdot |$ to $K(X)$ by
multiplicativity. Thus for any $L(X)\in K(X)$,
$L(X)=\frac{F(X)}{G(X)}$, with $F(X),G(X)\in K[X]$, $G(X)\neq 0$,
let $|L(X)|=\frac{|F(X)|}{|G(X)|}$. Let us remark that for any
non-zero element $u$ of $K[X]$ one has $|u|\geq 1$. In particular,
$R(\overline{g}^{m},F_{1})$ being a non-zero element of $K[X]$, we
have
\begin{equation}
|R(\overline{g}^{m},F_{1})|\geq 1.\label{ec4}
\end{equation} Next,
we estimate $|R(\overline{g}^{m},F_{1})|$ in a different way. Let
$ \overline{K(X)}$ be a fixed algebraic closure of $K(X)$, and let
us fix an extension of the absolute value $|\cdot |$ to
$\overline{K(X)}$, which we will also denote by $|\cdot |$.
Consider now the factorizations of $ \overline{g}(X,Y)$,
$\overline{g}^{m}(X,Y)$ and $F_{1}(X,Y)$ over $\overline{K(X)}$.
Say
\begin{align*}
\overline{g}(X,Y) &= \overline{b}_{n}(Y-\xi _{1})\cdots (Y-\xi _{n}), \\[.3pc]
\overline{g}^{m}(X,Y) &= \overline{b}_{n}^{m}(Y-\xi _{1})^{m}\cdots (Y-\xi
_{n})^{m}\\
\intertext{and}
F_{1}(X,Y) &= t_{1}(Y-\theta _{1})\cdots (Y-\theta _{r}),
\end{align*}
with $\xi _{1},\ldots ,\xi _{n},\theta _{1},\ldots ,\theta
_{r}\in $ $ \overline{K(X)}$. Here $1\leq r\leq mn-1$, by
our assumption that $\deg _{Y}F_{1}(X,Y)\geq 1$ and $\deg
_{Y}F_{2}(X,Y)\geq 1$. Then
\begin{equation}
|R(\overline{g}^{m},F_{1})|=\left|\overline{b}_{n}^{mr}t_{1}^{mn}\prod_{1\leq
i\leq n}\prod_{1\leq j\leq r}(\xi _{i}-\theta
_{j})^{m}\right|=|t_{1}|^{mn}\prod_{1\leq j\leq
r}|\overline{g}^{m}(X,\theta _{j})|.  \label{ec5}
\end{equation}
The fact that $F_{1}(X,\theta _{j})=0$ for any $j\in \{1,\ldots
,r\}$ also implies that $f(X,g(X,\theta _{j}))=0$, and so
\begin{align*}
h(X,\theta _{j}) &= f(X,g(X,\theta _{j}))-a_{m}(X)g^{m}(X,\theta
_{j}) \\[.3pc]
&= -a_{m}(X)b^{m}(X)\overline{g}^{m}(X,\theta _{j}).
\end{align*}
Since $b(X)$ is a non-zero element of $K[X]$, one has $|b(X)|\geq
1$. We deduce that
\begin{equation}
|\overline{g}^{m}(X,\theta _{j})|=\frac{|h(X,\theta _{j})|}{
|a_{m}(X)b^{m}(X)|}\leq \frac{|h(X,\theta _{j})|}{|a_{m}(X)|}.
\label{ec6}
\end{equation}
By combining (\ref{ec5}) and (\ref{ec6}) we find that
\begin{equation}
|R(\overline{g}^{m},F_{1})|\leq
\frac{|t_{1}|^{mn}}{|a_{m}|^{r}}\prod_{1\leq j\leq r}|h(X,\theta
_{j})|.  \label{ec7}
\end{equation}
We now proceed to find an upper bound for $|h(X,\theta
_{j})|$, for $1\leq j\leq r$. In order to do this, we
first use the identity
\begin{equation*}
h(X,Y)=a_{0}(X)+a_{1}(X)g(X,Y)+\cdots +a_{m-1}(X)g(X,Y)^{m-1}
\end{equation*}
to obtain
\begin{align}
|h(X,\theta _{j})| &= |a_{0}(X)+a_{1}(X)g(X,\theta _{j})+\cdots
+a_{m-1}(X)g(X,\theta _{j})^{m-1}|\nonumber \\[.1pc]
&\leq \max_{0\leq k\leq m-1}|a_{k}(X)|\cdot |g(X,\theta
_{j})|^{k}, \label{ec8}
\end{align}
for $1\leq j\leq r$. Next, we consider the factorization of
$f(X,Y)$ over $ \overline{K(X)}$, say
\begin{equation*}
f(X,Y)=a_{m}(X)(Y-\lambda _{1})\cdots (Y-\lambda _{m}),
\end{equation*}
with $\lambda _{1},\ldots ,\lambda _{m}\in \overline{K(X)}$. For
any $i\in \{1,\ldots ,m\}$ one has
\begin{equation}
0=f(X,\lambda _{i})=a_{0}(X)+a_{1}(X)\lambda _{i}+\cdots
+a_{m}(X)\lambda _{i}^{m}.  \label{ec9}
\end{equation}
By (\ref{ec9}) we see that
\begin{align}
|a_{m}(X)|\cdot |\lambda _{i}^{m}| &= |a_{0}(X)+a_{1}(X)\lambda
_{i}+\cdots + a_{m-1}(X)\lambda _{i}^{m-1}|  \nonumber \\[.1pc]
&\leq \max_{0\leq c\leq m-1}|a_{c}(X)|\cdot |\lambda _{i}|^{c}.
\label{ec10}
\end{align}
For any $i\in \{1,\ldots ,m\}$ let us select an index $c_{i}\in
\{0,\ldots ,m-1\}$ for which the maximum is attained on the right
side of (\ref{ec10}). We then have $|a_{m}(X)|\cdot |\lambda
_{i}|^{m}\leq |a_{c_{i}}(X)|\cdot |\lambda _{i}|^{c_{i}}$, and so
\begin{equation}
|\lambda _{i}|\leq \left( \frac{|a_{c_{i}}(X)|}{|a_{m}(X)|}\right)
^{1/(m-c_{i})}.  \label{ec11}
\end{equation}
We now return to (\ref{ec8}). Fix a $j\in \{1,\ldots ,r\}$. In
order to provide an upper bound for $|h(X,\theta _{j})|$, it is
sufficient to find an upper bound for $|g(X,\theta _{j})|$. Recall
that $f(X,g(X,\theta _{j}))=0$. Therefore there exists an $i\in
\{1,\ldots ,m\}$, depending on $j$, for which $g(X,\theta
_{j})=\lambda _{i}$. Then, by (\ref{ec11}) we obtain
\begin{equation}
|g(X,\theta _{j})|\leq \left(
\frac{|a_{c_{i}}(X)|}{|a_{m}(X)|}\right) ^{1/(m-c_{i})}\leq
\max_{1\leq v\leq m}\left( \frac{|a_{m-v}(X)|}{|a_{m}(X)|} \right)
^{1/v}.  \label{ec12}
\end{equation}
Inserting (\ref{ec12}) in (\ref{ec8}) we conclude that, uniformly
for $1\leq j\leq r$, one has
\begin{equation}
|h(X,\theta _{j})|\leq \max_{_{\substack{1\leq v\leq m \\ 0\leq
k\leq m-1}} }|a_{k}|\left( \frac{|a_{m-v}|}{|a_{m}|}\right)
^{k/v}. \label{ec13}
\end{equation}
Combining (\ref{ec13}) with (\ref{ec7}) we derive the inequality
\begin{equation*}
|R(\overline{g}^{m},F_{1})|\leq
\frac{|t_{1}|^{mn}}{|a_{m}|^{r}}\max_{_{\substack{1\leq v\leq m
\\ 0\leq k\leq m-1}}}\frac{|a_{k}|^{r}\cdot
|a_{m-v}|^{rk/v}}{|a_{m}|^{rk/v}},
\end{equation*}
which may be written as
\begin{equation}
|R(\overline{g}^{m},F_{1})|\leq |t_{1}|^{mn} \left(
\max_{_{\substack{1\leq v\leq m \\ 0\leq k\leq m-1}}}
\frac{|a_{k}|\cdot |a_{m-v}|^{k/v}}{ |a_{m}|^{1+k/v}}\right) ^{r}.
\label{ecu14}
\end{equation}
In what follows we are going to prove that
\begin{equation}
|t_{1}|^{mn}\max_{_{\substack{1\leq v\leq m \\ 0\leq k\leq
m-1}}}\frac{ |a_{k}|\cdot |a_{m-v}|^{k/v}}{|a_{m}|^{1+k/v}}<1,
\label{ecu15}
\end{equation}
which by (\ref{ecu14}) will contradict (\ref{ec4}), since $r\geq
1$. Using the definition of the absolute value $|\cdot |$, we
write the inequality (\ref{ecu15}) in the form
\begin{equation*}
\max_{_{\substack{1\leq v\leq m \\ 0\leq k\leq m-1}}}\rho
^{(1+\frac{k}{v} )\deg a_{m}-\deg a_{k}-\frac{k}{v}\deg
a_{m-v}}<\rho ^{mn\deg _{X}t_{1}},
\end{equation*}
which is equivalent to
\begin{equation}
\min_{_{\substack{1\leq v\leq m \\ 0\leq k\leq m-1}}}\left\{
\left(1+\frac{k}{v} \right)\deg a_{m}-\deg a_{k}-\frac{k}{v}\deg
a_{m-v}\right\} >mn\deg t_{1}.  \label{ec17}
\end{equation}
By combining (\ref{ec17}) with (\ref{ec2}), it will be sufficient
to prove that
\begin{align}
&\min_{_{\substack{1\leq v\leq m \\ 0\leq k\leq m-1}}}\left\{
\left(1+\frac{k}{v}\right)
\deg a_{m}-\deg a_{k}-\frac{k}{v}\deg a_{m-v}\right\} \nonumber
\\[.3pc]
&\quad\, >mn\deg d_{1}+m^{2}n\deg d_{2},\nonumber
\end{align}
or equivalently,
\begin{align}
&\min_{_{\substack{1\leq v\leq m \\ 0\leq k\leq
m-1}}}\left(1+\frac{k}{v}\right)\cdot \left\{ \deg
a_{m}-\frac{\deg a_{k}+\frac{k}{v}\deg a_{m-v}}{1+
\frac{k}{v}}\right\}   \nonumber \\[.3pc]
&\quad\, >mn\deg d_{1}+m^{2}n\deg d_{2}. \label{ecu18}
\end{align}
By our assumption on the size of $\deg a_{m}$ we have
\begin{equation*}
\deg a_{m}-H_{1}(f)>mn\deg d_{1}+m^{2}n\deg d_{2},
\end{equation*}
from which (\ref{ecu18}) follows, since
\begin{equation}
\max_{_{\substack{1\leq v\leq m \\ 0\leq k\leq m-1}}}\frac{\deg
a_{k}+ \frac{k}{v}\deg a_{m-v}}{1+\frac{k}{v}}\leq H_{1}(f).
\label{ecumax}
\end{equation}
This completes the proof of the first part of the theorem.
Assuming now that $f$ is irreducible over $K(X)$, the proof goes
as in the first part, except that now we have $\deg
_{Y}F_{1}=r\geq m$, since by Capelli's Theorem \cite{SCH}, the
degree in $Y$ of every irreducible factor of $f(X,g(X,Y))$ must be
a multiple of $m$. Therefore, instead of (\ref{ecu15}) one has to
prove that
\begin{equation*}
|t_{1}|^{n}\max_{_{\substack{ 1\leq v\leq m \\ 0\leq k\leq
m-1}}}\frac{ |a_{k}|\cdot |a_{m-v}|^{k/v}}{|a_{m}|^{1+k/v}}<1,
\end{equation*}
which is equivalent to
\begin{equation*}
\min_{_{\substack{ 1\leq v\leq m \\ 0\leq k\leq m-1}}}\left\{
\left(1+\frac{k}{v} \right)\deg a_{m}-\deg a_{k}-\frac{k}{v}\deg
a_{m-v}\right\} >n\deg t_{1}.
\end{equation*}
By combining this inequality with (\ref{ec2}), it will be
sufficient to prove that
\begin{align}
&\min_{_{\substack{ 1\leq v\leq m \\ 0\leq k\leq
m-1}}}\left(1+\frac{k}{v}\right)\cdot \left\{ \deg
a_{m}-\frac{\deg a_{k}+\frac{k}{v}\deg a_{m-v}}{1+\frac{k}{v}}
\right\}   \nonumber \\[.3pc]
&\quad\, >n\deg d_{1}+mn\deg d_{2}. \label{ecu19}
\end{align}
Finally, our assumption that $\deg a_{m}-H_{1}(f)>n\deg
d_{1}+mn\deg d_{2}$ together with (\ref{ecumax}) imply
(\ref{ecu19}), which completes the proof of the theorem. \hfill
$\Box$\vspace{.7pc}

We end by noting that in the statement of Theorem~1, the
assumption on the size of $\deg a_{m}$, and the bound $\Omega
(a_{m}/d_{1})+m\Omega (b_{n}/d_{2})$ exhibited for the number of
factors do not depend on the first $n$ coefficients of $g$. So
these bounds remain the same once we fix $n,b_{n}(X)$ and
$d_{2}(X)$, and let $b_{0}(X),\ldots ,b_{n-1}(X)$ vary
independently.

\section*{Acknowledgements}

This work was partially supported by the CERES Program 3-28/2003
and the CNCSIS Grant 40223/2003 code 1603 of the Romanian Ministry
of Education and Research.

\end{document}